\documentclass{amsart}

\usepackage{amsmath,amsthm,amsfonts,amscd,amssymb}
\usepackage[pdfborder={0 0 0}]{hyperref}
\usepackage{enumerate,verbatim}
\usepackage{enumitem}
\usepackage[usenames,dvipsnames,svgnames,table]{xcolor}

\usepackage{color}
\usepackage{graphicx}
\usepackage{fullpage}
\usepackage{enumerate}
\DeclareMathOperator*{\argmin}{arg min}

\newcommand{\NN}{\mathbb{N}}
\newcommand{\ZZ}{\mathbb{Z}}

\newtheorem{thm}{Theorem}[section]
\newtheorem{cor}[thm]{Corollary}
\newtheorem{lem}[thm]{Lemma}
\newtheorem{prop}[thm]{Proposition}

\newtheorem{prob}{Problem}

\theoremstyle{definition}

\theoremstyle{remark}
\newtheorem{rem}{Remark}

\theoremstyle{remark}
\newtheorem{example}{Example}

\begin{document}

\def\A{{\mathcal A}}
\def\AA{{\mathfrak A}}
\def\B{{\mathcal B}}
\def\C{{\mathcal C}}
\def\D{{\mathcal D}}
\def\EE{{\mathfrak E}}
\def\F{{\mathcal F}}
\def\G{{\mathcal G}}
\def\x{{\mathcal H}}
\def\I{{\mathcal I}}
\def\II{{\mathfrak I}}
\def\J{{\mathcal J}}
\def\K{{\mathcal K}}
\def\kk{{\mathfrak K}}
\def\L{{\mathcal L}}
\def\LL{{\mathfrak L}}
\def\M{{\mathcal M}}
\def\mm{{\mathfrak m}}
\def\MM{{\mathfrak M}}
\def\N{{\mathcal N}}
\def\O{{\mathcal O}}
\def\OO{{\mathfrak O}}
\def\PP{{\mathfrak P}}
\def\R{{\mathcal R}}
\def\W{{\mathcal W}}
\def\PNR{{\mathcal P_N(\real)}}
\def\PMNR{{\mathcal P^M_N(\real)}}
\def\PdNR{{\mathcal P^d_N(\real)}}
\def\s{{\mathcal S}}
\def\V{{\mathcal V}}
\def\X{{\mathcal X}}
\def\Y{{\mathcal Y}}
\def\Z{{\mathcal Z}}
\def\H{{\mathcal H}}
\def\cee{{\mathbb C}}
\def\Nn{{\mathbb N}}
\def\pee{{\mathbb P}}
\def\que{{\mathbb Q}}
\def\QQ{{\mathbb Q}}
\def\real{{\mathbb R}}
\def\RR{{\mathbb R}}
\def\zed{{\mathbb Z}}
\def\ZZ{{\mathbb Z}}
\def\aaa{{\mathbb A}}
\def\ff{{\mathbb F}}
\def\II{{\mathbb I}}
\def\HDelta{{\it \Delta}}
\def\kk{{\mathfrak K}}
\def\qbar{{\overline{\mathbb Q}}}
\def\kbar{{\overline{K}}}
\def\ybar{{\overline{Y}}}
\def\kkbar{{\overline{\mathfrak K}}}
\def\ubar{{\overline{U}}}
\def\eps{{\varepsilon}}
\def\ahat{{\hat \alpha}}
\def\bhat{{\hat \beta}}
\def\k{{\nu}}
\def\gt{{\tilde \gamma}}
\def\h{{\tfrac12}}
\def\baa{{\boldsymbol \alpha}}
\def\be{{\boldsymbol e}}
\def\bei{{\boldsymbol e_i}}
\def\bc{{\boldsymbol c}}
\def\bdt{{\boldsymbol \delta}}
\def\bff{{\boldsymbol f}}
\def\bm{{\boldsymbol m}}
\def\bk{{\boldsymbol k}}
\def\by{{\boldsymbol y}}
\def\bi{{\boldsymbol i}}
\def\bl{{\boldsymbol l}}
\def\bq{{\boldsymbol q}}
\def\bu{{\boldsymbol u}}
\def\bt{{\boldsymbol t}}
\def\bs{{\boldsymbol s}}
\def\bv{{\boldsymbol v}}
\def\bw{{\boldsymbol w}}
\def\bx{{\boldsymbol x}}
\def\bbx{{\overline{\boldsymbol x}}}
\def\bX{{\boldsymbol X}}
\def\bz{{\boldsymbol z}}
\def\bwy{{\boldsymbol y}}
\def\bY{{\boldsymbol Y}}
\def\bL{{\boldsymbol L}}
\def\ba{{\boldsymbol a}}
\def\bb{{\boldsymbol b}}
\def\bth{{\boldsymbol \theta}}
\def\bet{{\boldsymbol\eta}}
\def\bxi{{\boldsymbol\xi}}
\def\bo{{\boldsymbol 0}}
\def\bone{{\boldsymbol 1}}
\def\bol{{\boldsymbol 1}_L}
\def\ep{\varepsilon}
\def\p{\boldsymbol\varphi}
\def\q{\boldsymbol\psi}
\def\rank{\operatorname{rank}}
\def\aut{\operatorname{Aut}}
\def\lcm{\operatorname{lcm}}
\def\sgn{\operatorname{sgn}}
\def\spn{\operatorname{span}}
\def\md{\operatorname{mod}}
\def\Norm{\operatorname{Norm}}
\def\dim{\operatorname{dim}}
\def\det{\operatorname{det}}
\def\Vol{\operatorname{Vol}}
\def\rk{\operatorname{rk}}
\def\ord{\operatorname{ord}}
\def\ker{\operatorname{ker}}
\def\div{\operatorname{div}}
\def\Gal{\operatorname{Gal}}
\def\GL{\operatorname{GL}}
\def\SNR{\operatorname{SNR}}
\def\WR{\operatorname{WR}}
\def\IWR{\operatorname{IWR}}
\def\scg{\operatorname{\left< \Gamma \right>}}
\def\swrh{\operatorname{Sim_{WR}(\Lambda_h)}}
\def\ch{\operatorname{C_h}}
\def\cht{\operatorname{C_h(\theta)}}
\def\scgt{\operatorname{\left< \Gamma_{\theta} \right>}}
\def\scgmn{\operatorname{\left< \Gamma_{m,n} \right>}}
\def\gat{\operatorname{\Omega_{\theta}}}
\def\Obar{\operatorname{\overline{\Omega}}}
\def\Lbar{\operatorname{\overline{\Lambda}}}
\def\mn{\operatorname{mn}}
\def\disc{\operatorname{disc}}
\def\rot{\operatorname{rot}}
\def\Prob{\operatorname{Prob}}
\def\co{\operatorname{co}}
\def\ot{\operatorname{o_{\tau}}}
\def\Aut{\operatorname{Aut}}
\def\Mat{\operatorname{Mat}}
\def\SL{\operatorname{SL}}
\def\id{\operatorname{id}}
\def\supp{\operatorname{supp}}

\title{An algebraic perspective on integer sparse recovery}
\author{Lenny Fukshansky}\thanks{Fukshansky acknowledges support of the Simons Foundation grant \#519058}
\author{Deanna Needell}\thanks{Needell acknowledges support of NSF CAREER DMS $\#$1348721, NSF BIGDATA DMS $\#$1740325, and MSRI NSF DMS $\#$1440140}
\author{Benny Sudakov}\thanks{Sudakov acknowledges support of SNSF grant 200021-175573}

\address{Department of Mathematics, 850 Columbia Avenue, Claremont McKenna College, Claremont, CA 91711, USA}
\email{lenny@cmc.edu}
\address{Department of Mathematics, University of California at Los Angeles, 520 Portola Plaza, Los Angeles, CA 90095, USA}
\email{deanna@math.ucla.edu}
\address{Department of Mathematics, ETH Z\"urich, HG G 65.1 Ramistrasse 101, 8092 Z\"urich, Switzerland}
\email{benjamin.sudakov@math.ethz.ch}


\begin{abstract}
Compressed sensing is a relatively new mathematical paradigm that shows a small number of linear measurements are enough to efficiently reconstruct a large dimensional signal under the assumption the signal is \textit{sparse}. Applications for this technology are ubiquitous, ranging from wireless communications to medical imaging, and there is now a solid foundation of mathematical theory and algorithms to robustly and efficiently reconstruct such signals. However, in many of these applications, the signals of interest do not only have a sparse representation, but have other structure such as lattice-valued coefficients. While there has been a small amount of work in this setting, it is still not very well understood how such extra information can be utilized during sampling and reconstruction. Here, we explore the problem of integer sparse reconstruction, lending insight into when this knowledge can be useful, and what types of sampling designs lead to robust reconstruction guarantees. We use a combination of combinatorial, probabilistic and number-theoretic methods to discuss existence and some constructions of such sensing matrices with concrete examples. We also prove sparse versions of Minkowski's Convex Body and Linear Forms theorems that exhibit some limitations of this framework.
\end{abstract}
\maketitle

\section{Introduction}
\label{intro}
Initially motivated by a seemingly wasteful signal acquisition paradigm, compressed sensing has become a broad body of scientific work spanning across the disciplines of mathematics, computer science, statistics, and electrical engineering \cite{RefWorks:45,RefWorks:373}. Described succinctly, the main goal of compressed sensing is sparse recovery -- the robust \textit{reconstruction} (or \textit{decoding}) of a sparse signal from a small number of linear measurements. That is, given a signal $\bx\in\real^d$, the goal is to accurately reconstruct $\bx$ from its noisy measurements 
\begin{equation}\label{noisy}
\bb = A\bx + \be\in\real^m.
\end{equation}
Here, $A$ is an underdetermined matrix $A\in\real^{m\times d}$ ($m\ll d$), and $\be\in\real^m$ is a vector modeling noise in the system. Since the system is highly underdetermined, it is ill-posed until one imposes additional constraints, such as the signal $\bx$ obeying a sparsity model. We say $\bx$ is $s$-sparse when it has at most $s$ nonzero entries: 
\begin{equation}\label{sparsity}
\|\bx\|_0 := |\supp(\bx)| = |\{i : x_i \ne 0\}| \leq s\ll d.
\end{equation}
Clearly, any matrix $A$ that is one-to-one on $2s$-sparse signals will allow reconstruction in the noiseless case ($\be=0$). However, compressed sensing seeks the ability to reconstruct \textit{efficiently} and \textit{robustly}; one needs a computationally feasible reconstruction method, and one that allows accurate reconstruction even in the presence of noise. Fortunately, there is now a large body of work showing such methods are possible even when $m$ is only logarithmic in the ambient dimension, $m\approx s\log(d)$ \cite{RefWorks:45,RefWorks:373}. Typical results rely on notions like incoherence, null-space property or the restricted isometry property \cite{RefWorks:48}, which are quantitative properties of the matrix $A$ slightly stronger than simple injectivity. Under such assumptions, greedy (e.g. \cite{Paper9,NeedeT_CoSaMP,PaperIHT,RefWorks:150})  and optimization-based (e.g. \cite{RefWorks:48,RefWorks:47}) approaches have been designed and analyzed that efficiently produce an estimation $\hat{\bx}$ to an $s$-sparse signal $\bx\in\real^d$ from its measurements $\bb = A\bx + \be\in\real^m$ that satisfies
\begin{equation}\label{robust}
\|\hat{\bx} - \bx\| \lesssim \|\be\|,
\end{equation}
where $\lesssim$ hides only constant factors and $\|\cdot\|$ will always denote the Euclidean norm. Although this body of work has blossomed into many other directions based on practical motivations, there is very little understanding about the role of lattice-valued signals in this paradigm. This is especially troubling given the abundance of applications in which the signal is known to have lattice-valued coefficients, such as in wireless communications \cite{rossi2014spatial}, collaborative filtering \cite{davenport2016overview}, error correcting codes \cite{candes2005error}, and many more.  Initial progress in this setting includes results for dense (not sparse) $\pm1$ signals \cite{mangasarian2011probability}, binary sparse signals \cite{RefWorks:288,stojnic2010recovery}, finite-alphabet sparse signals \cite{tian2009detection,zhu2011exploiting}, and generalized lattice-valued signals \cite{flinth2017promp}.  The latter two categories are most relevant to our work; \cite{tian2009detection,zhu2011exploiting} propose modifications of the sphere decoder method that offer some empirical advantages but lack a rigorous theory. The recent work \cite{flinth2017promp} provides some theoretical guarantees for the greedy method OMP \cite{Paper9} initialized with a pre-processing step, and also shows that rounding the result given by $\ell_1$-minimization does not yield any improvements for many lattices. In this paper, our focus is not algorithmic but instead we aim to answer the questions (i) what kind of sensing matrices can be designed for lattice-valued sparse signals, and (ii) what are the limitations of the advantages one hopes to gain from knowledge that the signal is lattice-valued?  Our perspective in this work is thus \textit{algebraic}, and we leave algorithmic designs for such lattice-valued settings for future work.  We view our contribution as the foundation of an algebraic framework for lattice-valued signal reconstruction, highlighting both the potential and the limitations.
\bigskip

\section{Problem formulation and main results}
Let $m < d$ and first consider the noiseless consistent underdetermined linear system
\begin{equation}
\label{system1}
A\bx = \bb,
\end{equation}
where $A$ is an $m \times d$ real matrix and $\bb \in \real^d$. Let us first consider when this system has a unique solution~$\bx$. 
Notice that if $\bx$ and $\bwy$ are two different such solution vectors, then
$$A(\bx-\bwy) = \bo,$$
i.e. the difference vector $\bx-\bwy \in \NN(A)$, the null-space of the matrix $A$. 

Let us write $\ba_1,\dots,\ba_d \in \real^m$ for the column vectors of the matrix $A$. A vector $\bz \in \NN(A)$ if and only if
\begin{equation}
\label{system2}
A\bz = \sum_{i=1}^d z_i \ba_i = \bo,
\end{equation}
i.e. if and only if $\ba_1,\dots,\ba_d$ satisfy a linear relation with coefficients $z_1,\dots,z_n$. The uniqueness of solution to~\eqref{system1} (and hence our ability to decode the original signal) is equivalent to non-existence of nonzero solutions to~\eqref{system2}.

Since $d > m$, such solutions to~\eqref{system2} must exist. On the other hand, if we add some appropriate restriction on the solution vectors in question, then perhaps there will be no solutions satisfying this restriction. In other words, the idea is to ensure uniqueness of decoded signal by restricting the original signal space.  We can then formulate the following problem.

\begin{prob} \label{P1} Define a restricted $d$-dimensional signal space $X \subseteq \real^d$ and an $m \times d$ matrix $A$ with $m < d$ such that $A\bx \neq 0$ for any $\bx \in X$.
\end{prob}

A commonly used restriction is sparsity, defined in \eqref{sparsity}. Now, while~\eqref{system2} has nonzero solutions, it may not have any nonzero $s$-sparse solutions for sufficiently small~$s$. In addition to exploiting sparsity, one can also try taking advantage of another way of restricting the signal space. Specifically, instead of taking signals to lie over the field~$\real$, we can restrict the coordinates to a smaller subfield of~$\real$, for instance~$\que$. The idea here is that, while columns of our matrix $A$ are linearly dependent over $\real$, they may still be linearly independent over a smaller subfield. For instance, we have the following trivial observation.

\begin{lem} \label{trivial} Let $K \subsetneq \real$ be a proper subfield of the field of real numbers, and let $\alpha_1,\dots,\alpha_d \in \real$ be linearly independent over $K$. Define $A = ( \alpha_1\ \dots\ \alpha_d)$ be a $1 \times d$ matrix, then the equation $A \bx = 0$ has no solutions in $K^d$ except for $\bx=\bo$.
\end{lem}

Of course, when \eqref{system2} has no solutions, we can guarantee our system \eqref{system1} has a unique solution $\bx$ and can in theory will be able to decode successfully.  However, for practical concerns we want to be able to tolerate noise in the system and decode robustly as in \eqref{noisy}. Since in practice the noise $\be$ typically scales with the entries (or row or column norms) of $A$, we ask for the following two properties:
\renewcommand{\theenumi}{\roman{enumi}}
\begin{enumerate}
\item the entries of $A$ are uniformly bounded in absolute value
\item $\|A\bz\|$ is bounded away from zero for any vector $\bz\ne 0$ in our signal space (say, $\|A\bz\| > C$).
\end{enumerate}
For example, consider the signal space  
\begin{equation}
\label{Z2s}
\zed_{2s}^d := \left\{ \bz \in \zed^d : \|\bz\|_0 \leq 2s \right\},
\end{equation}
and suppose we wish to decode an $s$-sparse signal $\bx\in\zed_s^d$ from its noisy measurements $\bb = A\bx + \be$ where $\|\be\| \leq \frac{1}{2}C$. Suppose we decode (inefficiently) by selecting the signal $\by\in\zed_s^d$ minimizing $\|\bb - A\by\|$.  Then since $\bx\in\zed_s^d$ is such that $\|\bb - A\bx\| = \|\be\| \leq \frac{1}{2}C$, it must be that the decoded vector $\by$ satisfies $\|\bb - A\by\| \leq \frac{1}{2}C$ as well. Therefore, $\|A\by - A\bx\| \leq \|\bb - A\by\| + \|\bb - A\bx \| \leq C$. Then (noting that $\bx-\by\in\zed_{2s}^d$)  by (ii), this guarantees that $\by=\bx$ so our decoding was successful. Hence we consider the following optimization problem (which we will want to use with $s'=2s$).

\begin{prob} \label{P2}
Construct an $m \times d$ matrix $A$ with $m < d$ such that
$$|A| := \max \{ |a_{ij}| : 1 \leq i \leq m, 1 \leq j \leq d\} \leq C_1,$$
and for every nonzero $\bx \in \zed_{s'}^d$,
$$\|A\bx\| \geq C_2,$$
where $C_1,C_2 > 0$. 
\end{prob}

\noindent
In this paper, we discuss existence and construction of such matrices. Here is our first result.

\begin{thm} \label{main} There exist $m \times d$ integer matrices $A$ with $m < d$ and bounded $|A|$ such that for any nonzero $\bx \in \zed_s^d$, $0 < s \leq m$, 
\begin{equation}
\label{Ax-bnd}
\|A\bx\| \geq 1.
\end{equation}
In fact, for sufficiently large $m$, there exist such matrices with
\begin{equation}
\label{md-1}
|A|=1 \text{ and } d=1.2938\, m,
\end{equation}
and there also exist such matrices with
\begin{equation}
\label{md-2}
|A|=k \text{ and } d=\Omega(\sqrt{k}\, m).
\end{equation}
On the other hand, for any integers $m \geq 3$, $k \geq 1$ and $m \times d$ integer matrix $A$ with $|A|=k$ satisfying~\eqref{Ax-bnd} for all $s \leq m$, we must have
\begin{equation}
\label{m-bnd}
d \leq (2k^2+2) (m-1) + 1.
\end{equation}
\end{thm}

\begin{rem} Notice that in situations when one needs to have the bound~\eqref{Ax-bnd} replaced by a stronger bound $\|A\bx\| \geq \ell$ for some $\ell > 1$, this can be achieved by simply multiplying $A$ by $\ell$, of course at the expense of making~$|A|$ larger, but only by the constant factor of~$\ell$.
\end{rem}
\medskip

\noindent
We discuss the dependence between $m$, $d$ and $|A|$ in more detail and prove Theorems~\ref{main} in Section~\ref{integer}. In Section~\ref{new_mtrx}, we extend this matrix construction over number fields, proving the following corollary of Theorem~\ref{main}.

\begin{cor} \label{main_cor} Let $B$ be the $d \times m$-transpose of a matrix satisfying~\eqref{Ax-bnd} as guaranteed by Theorem~\ref{main}. Let $\theta$ be an algebraic integer of degree $m$, and let $\theta = \theta_1, \theta_2,\dots, \theta_m$ be its algebraic conjugates. For each $1 \leq i \leq m$, let $\bth_i = (1\ \theta_i\ \dots \theta_i^{m-1})^{\top}$, compute the $d \times m$ matrix
$$B \begin{pmatrix} \bth_1 & \dots & \bth_m \end{pmatrix},$$
and let $A$ be its transpose. Then $|A| = O \left( |B| m \right)$, for any $\bx \in \zed_s^d$, $0 < s \leq m$, $\|A\bx\| \geq \sqrt{m}$ and the vector $A\bx$ has all nonzero coordinates.
\end{cor}

\noindent
In Section~\ref{recon} we discuss an algorithm for reconstructing the original sparse signal~$\bx$ from its measurement $A\bx + \be$, where $\be$ is the error vector of Euclidean norm $< \sqrt{m}/2$. We show that the complexity of this algorithm is the same as that of the Closest Vector Problem (CVP) in~$\real^m$.

While these results show the existence of matrices $A$ such that $\|A\bx\|$ is bounded away from $\bo$ on sparse vectors, it is also clear that for any $m \times d$ matrix $A$ there exist sparse vectors with $\|A\bx\|$ not too large: for instance, if $\bx \in \zed^d$ is a standard basis vector, then
\begin{equation}
\label{naive}
\|A\bx\| \leq \sqrt{m}\ |A|.
\end{equation}
In Section~\ref{sparse} we prove a determinantal upper bound on $\|A\bx\|$ in the spirit of the Geometry of Numbers.

\begin{thm} \label{main2} Let $A$ be an $m \times d$ real matrix of rank $m \leq d$, and let $A'$ be the $d \times m$ real matrix so that $AA'$ is the $m \times m$ identity matrix. There exists a nonzero point $\bx \in \zed_m^d$ such that
\begin{equation}
\label{thm_2-4}
\|A\bx\| \leq \sqrt{m}\ \left| \det \left( (A')^{\top} A' \right) \right|^{-1/2m}.
\end{equation}
\end{thm}

\noindent
We prove this result by deriving sparse versions of Minkowski's Convex Body and Linear Forms Theorems for parallelepipeds. There are many situations when the bound of Theorem~\ref{main2} and the naive bound $\sqrt{m}\ |A|$ are comparable, but there are also many cases when the bound of~\eqref{thm_2-4} is substantially better than that of~\eqref{naive}. We demonstrate several such examples at the end of Section~\ref{sparse}. We are now ready to proceed.
\bigskip

\section{An integer matrix}
\label{integer}

Let us fix positive integers $s \leq m \leq d$. Let $\ba_1,\dots,\ba_d \in \zed^m$ be a collection of vectors such that no $m$ of them are linearly dependent. Define
\begin{equation}
\label{matrix_B}
A = \begin{pmatrix} \ba_1 & \dots & \ba_d \end{pmatrix},
\end{equation}
i.e., $A$ is an $m \times d$ integer matrix with column vectors $\ba_1,\dots,\ba_d$. Write $[d] := \{1,\dots,d\}$ and let $I \subset [d]$ be a subset of cardinality $m$. Let $A_I$ be an $m \times m$ submatrix of $A$, consisting of the columns indexed by elements of~$I$. The determinant $\det A_I$ is called the corresponding {\it Pl\"ucker coordinate} of $A$, and the set of all Pl\"ucker coordinates of $A$ is
\begin{equation}
\label{plucker}
\nu(A) = \{ \det A_I : I \subset [d], |I| = m \}.
\end{equation}
Notice that the condition that no $m$ of $\ba_1,\dots,\ba_d$ are linearly dependent is equivalent to the condition that all Pl\"ucker coordinates of $A$ are nonzero.

\begin{lem} \label{int_bnd} Let $s\leq m$. For every nonzero $\bx \in \zed_s^d$,
\begin{equation}
\label{Ax_norm}
\|A\bx\| \geq 1.
\end{equation}
\end{lem}

\proof
Let $\bx \in \zed_s^d$, then at most $s$ of coordinates of $\bx$ are nonzero, assume these are $x_{j_1},\dots,x_{j_s}$ for some $1 \leq j_1 < \dots < j_s \leq d$. Then
$$A\bx = x_{j_1}\ba_{j_1} + \dots + x_{j_s}\ba_{j_s} \neq \bo,$$
since $s \leq m$ and no $m$ vectors among $\ba_1,\dots,\ba_d$ can be linearly dependent. Since $A\bx$ is a nonzero integer vector, its Euclidean norm has to be at least 1.
\endproof
\smallskip
Since $m$ corresponds to the compressed dimension, we typically fix $d$ and ask for $m$ to be as small as possible, or equivalently, fix $m$ and ask how large $d$ can be.  It is clear that for any fixed $m$ we can take $d$ to be as large as possible, however this will force $|A|$ to grow. A simple argument shows that for any $m$ and $d$ it is always possible to construct an $m \times d$ integer matrix $A = (a_{ij})$ with all Pl\"ucker coordinates nonzero and
$$|A| \leq \frac{1}{2} \binom{d-1}{m-1}.$$
Indeed, we can let $P(A) = P(a_{ij})$ to be the product of determinants of all $m \times m$ submatrices, and then notice that this polynomial cannot vanish ``too much". Specifically, $P$ is a polynomial in $md$ variables $a_{ij}$ and its degree in
each of the variables is $\binom{d-1}{m-1}$, since each determinant of an $m \times m$ submatrix is linear in each of the
variables $a_{ij}$ and each column of $A$ is present in $t := \binom{d-1}{m-1}$ such submatrices. Since $P$ is not identically zero, it cannot vanish on all of $\{ -t/2,\dots,0,\dots,t/2 \}^{md}$ (see, for instance Lemma~2.1 of~\cite{alon}). Hence there must exist a matrix $A = (a_{ij})$ with all the entries $\leq t/2$ in absolute value such that $P(A) \neq 0$.

On the other hand, our Theorem~\ref{main} implies a much better bound on $|A|$ in terms of $m$ and $d$: notice that~\eqref{md-2} guarantees the existence of an $m \times d$ integer matrix $A$ with all nonzero Pl\"ucker coordinates so that $|A| = O(d/m)^2$. Let us now turn to the proof of this theorem.

\proof[Proof of Theorem~\ref{main}]
As we discussed above, condition~\eqref{Ax-bnd} is equivalent to saying that all Pl\"ucker coordinates of $A$ are nonzero. The fact that such $m \times d$ integer matrices exist with $m < d$ and $|A|=1$ is immediate: for any $m$ take $m \times (m+1)$ matrix with first $m$ columns being standard basis vectors in~$\real^m$ and the last column being $(1\ \dots\ 1)^{\top}$. To obtain better results stated in ~\eqref{md-1} and~\eqref{md-2} we use probabilistic arguments. 

To prove~\eqref{md-1} we need the following result from \cite{random} (Corollary 3.1).
Let $M_m$ be an $m \times m$ matrix whose entries are independent copies of a random variable $\mu$ taking value $0$ with probability $1/2$ and values $-1$ or $1$ with probability $1/4$.
Then the probability that matrix $M_m$ is singular is at most $(1/2-o(1))^m$. Form an $m \times d$ random matrix $A$ by taking its entries to be independent copies of $\mu$. Note that $|A|=1$ and any $m$ columns of $A$ form a matrix distributed according to $M_m$. Therefore the probability that any $m \times m$ submatrix of $A$ is singular is at most $(1/2-o(1))^m$.
Since the number of such submatrices is ${d \choose m}$ we have (by union bound) that the probability that $A$ contains an  $m \times m$ singular submatrix is at most 
$${d \choose m}(1/2-o(1))^m.$$
To bound this probability we use the following well known estimate on the binomial coefficients (see, e.g, Chapter~15 of~\cite{AS}). Let
$$H(p)=-p \log_2 p -(1-p) \log_2 (1-p)$$
be the binary entropy function, then ${d \choose pd} \leq 2^{H(p)d}$.
Using this estimate together with $d=1.2938m$, one can easily show that ${d \choose m}(1/2-o(1))^m < 1$. Thus with positive probability $A$ does not have singular $m \times m$ submatrices. This implies that there exists such a matrix $A$ so that for any $\bx \in \zed_s^d$, $0 < s \leq m$, $\|A\bx\| >0$; hence $\|A\bx\| \geq 1$ since it is an integer.

To prove~\eqref{md-2} we need another result from \cite{random} (Corollary 3.3). It says that if $N_m$ is an $m \times m$ matrix whose entries are independent copies of a random variable $\mu'$ taking uniformly one of the $2k+1$ integer values from $\{-k, \cdots, k\}$, then the probability that $N_m$ is singular is at most $(1/\sqrt{2k}-o(1))^m$.
Consider an $m \times d$ random matrix $A$ whose entries are independent copies of $\mu'$. Then $|A|=k$ and 
the probability that any $m \times m$ submatrix of $A$ is singular is at most $(1/\sqrt{2k}-o(1))^m$.
Since the number of such submatrices is ${d \choose m}$ we have that the probability that
$A$ contains an  $m \times m$ singular submatrix is at most 
$${d \choose m}(1/\sqrt{2k}-o(1))^m.$$
Using the estimate ${a \choose b} \leq (ea/b)^b$ for binomial coefficient and choosing $d$ to be a sufficiently small multiple of $\sqrt{k}m$ this probability can be made smaller than $1$. Thus with positive probability $A$ does not have singular $m \times m$ submatrices, implying again that for any $\bx \in \zed_s^d$, $0 < s \leq m$, $\|A\bx\| \geq 1$.

Finally, let us turn to the proof of~\eqref{m-bnd}. Let $m \geq 3$ and $k \geq 1$ be fixed integers. Let $A$ be an $m \times d$ integer matrix with $|A| \leq k$ such that all Pl\"ucker coordinates of $A$ are nonzero. We want to show that $d \leq (2k^2+2) (m-1) + 1$, or equivalently, $m\geq \frac{d-1}{2k^2+2}+1$. For any real vector $\bx$, we write $|\bx|$ for its sup-norm, i.e. maximum of absolute values of its coordinates. Let 
$$C_m(k) = \{ \bx \in \zed^m : |\bx| \leq k \},$$
then $|C_m(k)| = (2k+1)^m$. Let $\ell=(2k^2+2) (m-1) + 1$ and let $\bx_1,\dots,\bx_\ell$ be any $\ell$ vectors from $C_m(k)$. If there are $m$ vectors with the first or second coordinate equal to $0$, then they all lie in the same $(m-1)$-dimensional subspace. If not, there are $2k^2m$ vectors with the first two coordinates nonzero. Multiplying some of these vectors by $-1$, if necessary (does not change linear independence properties) we can assume that the all of them have positive first coordinate. Hence there are a total of $k \times 2k = 2k^2$ choices for the first two coordinates, so there must exist a subset of $m$ of these vectors that have these first two coordinates the same, let these be $\bx_1,\dots,\bx_m$. Then there exists a vector $\bwy = (a,b,0,\dots,0)^{\top} \in C_m(k)$ such that the vectors
$$\bz_1 = \bx_1-\bwy,\ \dots,\ \bz_m = \bx_m-\bwy$$
all have the first two coordinates equal to 0. This means that these vectors lie in an $(m-2)$-dimensional subspace
$$V = \{ \bz \in \real^m : z_1 = z_2 = 0 \}$$
of $\real^m$. Then let $V' = \spn_{\real} \{ V, \bwy \}$, so $\dim_{\real} V' = m-1$. On the other hand, $\bx_1,\dots,\bx_m \in V'$, and hence the $m \times m$ matrix with rows $\bx_1,\dots,\bx_m$ must have determinant equal to 0. Therefore in order for an $m \times d$ matrix $A$ with column vectors in $C_m(k)$ to have all nonzero Pl\"ucker coordinates, $d$ has to be no bigger than $(2k^2+2) (m-1) + 1$.
\endproof

\begin{rem} \label{CV}
Notice that our $m \times d$ matrix $A$ as in~\eqref{matrix_B} has to have the property that no $s$ of its column vectors $\ba_1,\dots,\ba_d \in \zed^m$ are linearly dependent, i.e. every $s$-dimensional subspace of~$\real^m$ contains at most~$s$ of these vectors. Then let $d$ be the maximum number of vectors in~$\zed^m$ of Euclidean norm~$\leq r$ such that every $s$-dimensional subspace of~$\real^m$ contains at most~$s$ of them. Corollary~7 of~\cite{cover} then states that
$$d \leq O_{m,s} \left( r^{\frac{m(m-s)}{m-1}} \right),$$
where $O_{m,s}$ means that the constant in the $O$-notation depends only on $m$ and $s$. Corollary~6 of this same paper suggests a similar in spirit, but somewhat weaker lower bound for~$d$. These results, however, focus on the dependence of the bounds on $r$, not on the parameters $m$ and $s$ which are of main interest to us.
\end{rem}

Equation~\eqref{md-1} of our Theorem~\ref{main} guarantees the existence of an $m \times d$ integer matrix satisfying~\eqref{Ax-bnd} with $d=1.2938\, m$ for sufficiently large $m$. It is an interesting question what is the optimal dependence of~$d$ on~$m$ and also whether one can construct explicitly matrices satisfying assertions in ~\eqref{md-1} and~\eqref{md-2}? For instance, we can construct a simple $3$-dimensional example with $|A|=1$ and $d=2m$.

\begin{example} \label{int_matrix} Let $m=3$, $d=6$, $k=1$, and define a $3 \times 6$ matrix
\begin{equation}
\label{int1}
A = \begin{pmatrix} 1 & 1 & 1 & 1 & 1 & 1 \\ 1 & 1 & 0 & 0 & -1 & -1 \\ 1 & 0 & 1 & -1 & 0 & -1 \end{pmatrix}.
\end{equation}
This matrix has $|A| = 1$ and one can easily check that all its Pl\"ucker coordinates are nonzero, as required. Then for $s \leq 3$ and any $\bx \in \zed_s^6$, $\|A\bx\| \geq 1$. Note also that the maximal Euclidean norm of the row vectors of $A$ is $\approx 2.45$ and its smallest singular value is $1/\sqrt{2}$.
\end{example}
\bigskip

\section{Algebraic matrix construction}
\label{new_mtrx}

Here we extend our construction over number fields, proving Corollary~\ref{main_cor}. Let $1 \leq m \leq d$ be rational integers, and let $K$ be a number field of degree $m$ over $\que$ with embeddings $\id = \sigma_1,\sigma_2,\dots,\sigma_m$, where by $\id$ we mean the identity map on~$K$. Write $\O_K$ for the ring of algebraic integers of $K$. Let $d > 1$ and $\alpha_1,\dots,\alpha_d \in \O_K$. Define 
\begin{equation}
\label{matrix_A}
A = \begin{pmatrix} \sigma_1(\alpha_1) & \dots & \sigma_1(\alpha_d) \\ \vdots & \dots & \vdots \\ \sigma_m(\alpha_1) & \dots & \sigma_m(\alpha_d) \end{pmatrix}
\end{equation}
to be an $m \times d$ matrix over $K$, and for each $1 \leq i \leq m$, let
\begin{equation}
\label{L_form}
L_i(x_1,\dots,x_d) = \sum_{j=1}^d \sigma_i(\alpha_j) x_j
\end{equation}
be the linear form with coefficients $\sigma_i(\alpha_1),\dots,\sigma_i(\alpha_d)$, corresponding to the $i$-th row of $A$.

\begin{lem} \label{norm_lemma} Let the notation be as above, and suppose that $\Z \subseteq \zed^d$ is a signal space such that $L_1(\bx) \neq 0$ for any $\bo \neq \bx \in \Z$. Then for each $\bx \in \Z$,
$$\left\| A \bx \right\| \geq \sqrt{m},$$
where $\|\ \|$ stands for the usual Euclidean norm.
\end{lem}

\proof
Notice that for a vector $\bx \in \real^d$,
$$A\bx = (L_1(\bx),\dots,L_m(\bx))^{\top}\ \Longrightarrow\ \left\| A \bx \right\| = \left( \sum_{i=1}^m |L_i(\bx)|^2 \right)^{1/2}.$$
For each $\bx \in \zed^d$, $L_i(\bx) = \sigma_i(L_1(\bx)) \in \O_K$, and hence the field norm $\NN_K$ of $L_1(\bx)$ is a rational integer, i.e. unless $L_1(\bx) = 0$,
\begin{equation}
\label{norm}
\left| \NN_K(L_1(\bx)) \right| = \prod_{i=1}^m |\sigma_i(L_1(\bx))| =  \prod_{i=1}^m | L_i(\bx)| \geq 1.
\end{equation}
Now suppose $\bo \neq \bx \in \Z$, then $L_1(\bx) \neq 0$, and hence $L_i(\bx) \neq 0$ for all $1 \leq i \leq m$. Then, combining~\eqref{norm} with AM-GM inequality, we obtain:
\begin{equation}
\label{nrm_bnd}
\frac{1}{m} \sum_{i=1}^m |L_i(\bx)|^2 \geq \left( \prod_{i=1}^m |L_i(\bx)|^2 \right)^{1/m} = \left| \NN_K(L_1(\bx)) \right|^{\frac{2}{m}} \geq 1.
\end{equation}
The result follows.
\endproof

With Lemma~\ref{norm_lemma} in mind, we can now propose the following explicit construction. Let  $\alpha_1,\dots,\alpha_d \in \O_K$ be such that no $m$ of them are linearly dependent over~$\que$. For this choice of $\alpha_i$'s, let $A$ be as in~\eqref{matrix_A} and $L_i$'s as in~\eqref{L_form}. Let $1 \leq s \leq m$, then for any $\bo \neq \bx \in \zed_s^d$, $L_1(\bx) \neq 0$, and hence Lemma~\ref{norm_lemma} implies that
\begin{equation}
\label{Am}
\left\| A \bx \right\| \geq \sqrt{m}.
\end{equation}

We now want to find specific constructions of such a matrix $A$ so that the absolute values of its entries are small. Let us start with a small basis for $K$; we can, for instance take a power basis, i.e. if $K=\que(\theta)$ for an algebraic integer $\theta$, then
$$1, \theta, \theta^2, \dots, \theta^{m-1}$$
is a basis for $K$ over $\que$. Write $\bth$ for the column vector $(1, \theta, \theta^2, \dots, \theta^{m-1})^{\top}$. Let $k$ be a positive integer and let $B$ be a $d \times m$ matrix with integer entries in the interval $[-k,k]$ such that all Pl\"ucker coordinates of $B$ are nonzero; in other words, $B$ is the transpose of a matrix of the type we constructed in Section~\ref{integer}. Then define
$$\baa = (\alpha_1, \dots, \alpha_d) = B \bth,$$
and with these $\alpha_1,\dots,\alpha_d$ define the matrix $A$ as in~\eqref{matrix_A}. Notice that $A$ is an $m \times d$ matrix that has precisely the property~\eqref{Am} we need, since $s \leq m$. Further, 
$$|A| \leq m |B| |\theta|^{m-1} \leq m k |\theta|^{m-1}.$$
This finishes the construction of matrices $A$ as in Corollary~\ref{main_cor}, hence proving this corollary.
\smallskip

\begin{example}\label{ex4} Let $m=3$, $d=6$, and take $K = \que(\theta)$, where $\theta=2^{1/3}$, then
$$\bth = \begin{pmatrix} 1 \\ \theta \\ \theta^2 \end{pmatrix}.$$
Let $k=1$ and take $B$ to be the transpose of the matrix~\eqref{int1} from Example~\ref{int_matrix}, i.e.
$$B = \begin{pmatrix} 1 & 1 & 1 \\ 1 & 1 & 0 \\ 1 & 0 & 1 \\ 1 & 0 & -1 \\ 1 & -1 & 0 \\ 1 & -1 & -1 \end{pmatrix}.$$
Define
$$\baa = B\bth = \begin{pmatrix} 1 + \theta + \theta^2 & 1 + \theta & 1+ \theta^2 & 1 - \theta^2 & 1 - \theta & 1 - \theta - \theta^2 \end{pmatrix}.$$
The number field $K$ has three embeddings, given by $\theta \mapsto \theta$, $\theta \mapsto \xi \theta$, and $\theta \mapsto \xi^2 \theta$, where $\xi = e^{\frac{2\pi i}{3}}$ is a third root of unity, i.e. $\theta$ is mapped to roots of its minimal polynomial by injective field homomorphisms that fix~$\que$. Hence we get the following $3 \times 6$ matrix, as in~\eqref{matrix_A}:
$$A = \begin{pmatrix} 1 + \theta + \theta^2 & 1 + \theta & 1+ \theta^2 & 1 - \theta^2 & 1 - \theta & 1 - \theta - \theta^2 \\
1 + \xi \theta + \xi^2 \theta^2 & 1 + \xi \theta & 1+ \xi^2 \theta^2 & 1 - \xi^2 \theta^2 & 1 - \xi \theta & 1 - \xi \theta - \xi^2 \theta^2 \\
1 + \xi^2 \theta + \xi \theta^2 & 1 + \xi^2 \theta & 1+ \xi \theta^2 & 1 - \xi \theta^2 & 1 - \xi^2 \theta & 1 - \xi^2 \theta - \xi \theta^2
\end{pmatrix}$$
with $|A| \leq 3 \sqrt[3]{2}$ and $\left\| A \bx \right\| \geq \sqrt{3}$ for every $\bx \in \zed_s^6$, $s \leq 3$. 

\begin{rem}Note that every 3-column submatrix of this matrix $A$ is full-rank, and so the matrix contains no $3$-sparse signal (integer or real-valued) in its null-space. However, the smallest singular value of this matrix $A$ is $0.2736 \ll \sqrt{3}$; thus, one can only guarantee the norm bound for integer valued signals. Figure \ref{fig2} shows results for classical methods and their reconstruction performance using this matrix. These methods include L1-minimization \cite{RefWorks:47}, Orthogonal Matching Pursuit (OMP) \cite{Paper9}, simple least-squares (LS solves min$_{\bz}$ $\|A\bz-\bb\|$), and simple Hard Thresholding (which estimates the support as the largest [in magnitude] entries of $A^T\bb$ and then applies least-squares using that submatrix). These plots highlight the fact that although theoretical reconstruction is guaranteed, novel (efficient) reconstruction methods need to be created for these types of matrices and signals.
\end{rem}

\begin{rem}\label{workitout}
Let us remark on how we would use Lemma \ref{norm_lemma} in practice to guarantee robust recovery of sparse signals.  We will continue the illustration with Example \ref{ex4} for concreteness. Note that this example constructs a $3\times 6$ matrix $A$ with bounded entries such that for any non-zero $\by \in\zed_3^6$ we guarantee that $\left\| A \by \right\| \geq \sqrt{3}$. So, let $\bx\in\zed_1^6$ be given, and take noisy measurements $\bb = A\bx + \be$ where the noise obeys $\|\be\| < \frac{\sqrt{3}}{2}\approx 0.866$. Now, suppose we attempt to recover $\bx$ by (again, inefficiently) selecting the $1$-sparse integer vector $\hat{\bx}$ that minimizes $\|\bb-A\hat{\bx}\|$:
$$
\hat{\bx} := \argmin_{\by\in\zed_1^6} \|\bb-A\by\|.
$$
Then we have $\bx-\hat{\bx}\in\zed_2^6\subset\zed_3^6$ and $\|A\bx-A\hat{\bx}\| \leq \|\bb-A\bx\| + \|\bb-A\hat{\bx}\| < \frac{\sqrt{3}}{2} + \frac{\sqrt{3}}{2} = \sqrt{3}$. But by the construction of $A$ this must mean that $\hat{\bx} = \bx$, so we have reconstructed $\bx$ exactly.  See Figure \ref{therealfig2} for reconstruction results (using the inefficient CVP method), that shows we seem to tolerate noise slightly above this value of $\sqrt{3}/2$.
\end{rem}


\begin{figure}
\includegraphics[width=3in]{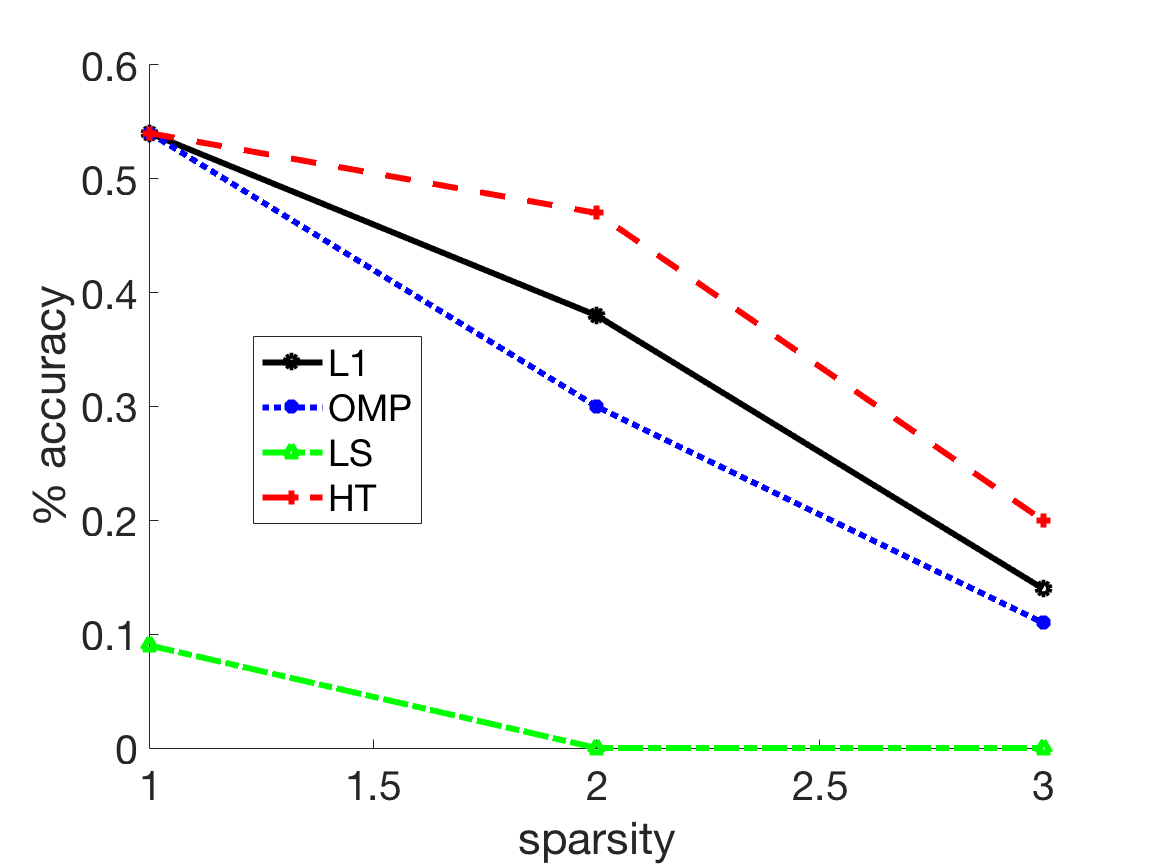}\quad\includegraphics[width=3in]{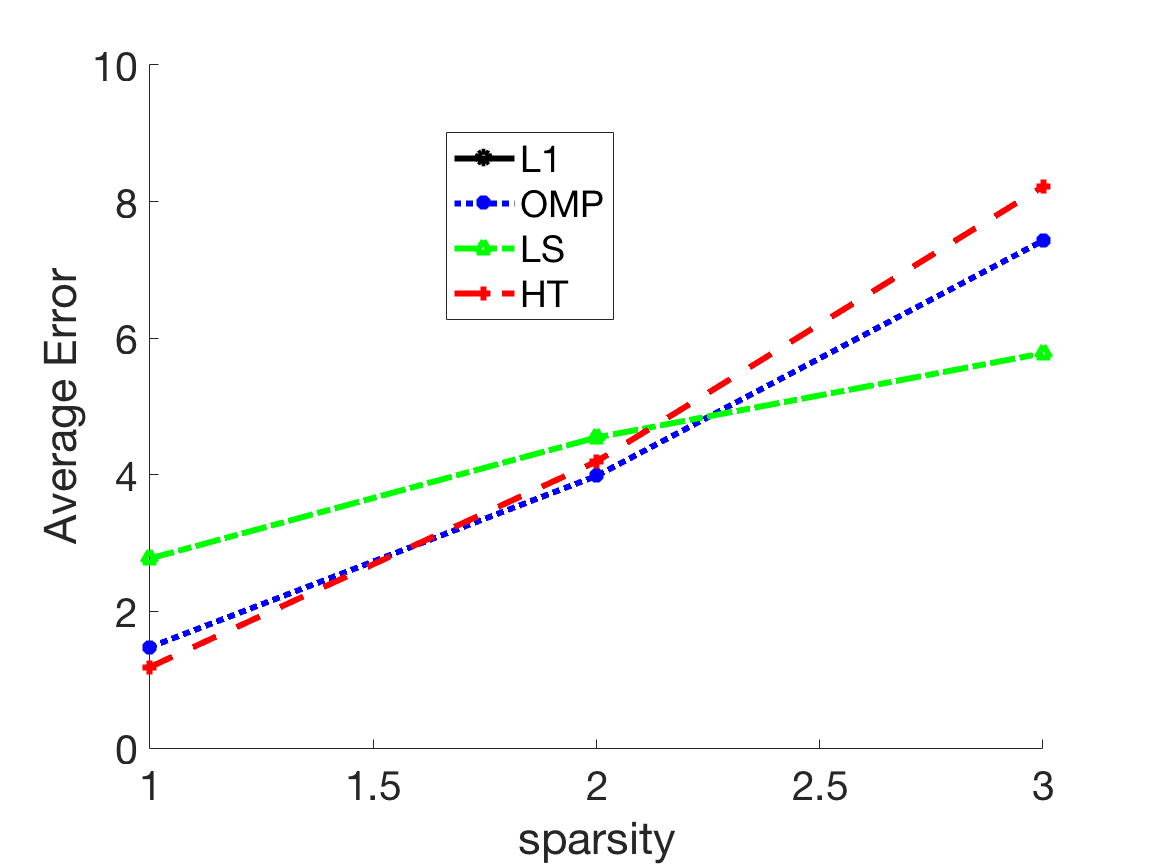}
\caption{Reconstruction results for the matrix described in Example \ref{ex4} using L1-minimization (L1), Orthogonal Matching Pursuit (OMP), Hard Thresholding (HT) and least-squares (LS). Signals were generated with random support, and entries followed a Gaussian distribution with mean 0 and variance $\sigma^2=25$ and were then rounded to the nearest integer. Right: Average (L2) reconstruction error when zero-mean Gaussian noise was added to the measurements $y=Ax$ before reconstruction, with $\sigma=0.5$ (L1 results not shown due to lack of convergence). }\label{fig2}
\end{figure}

\begin{figure}
\includegraphics[width=3in]{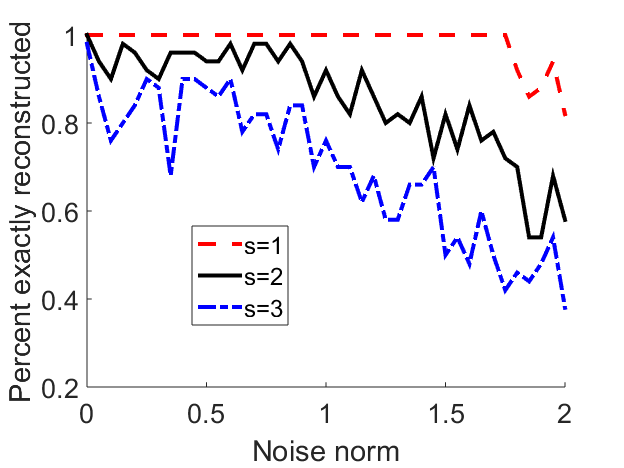}\quad\includegraphics[width=3in]{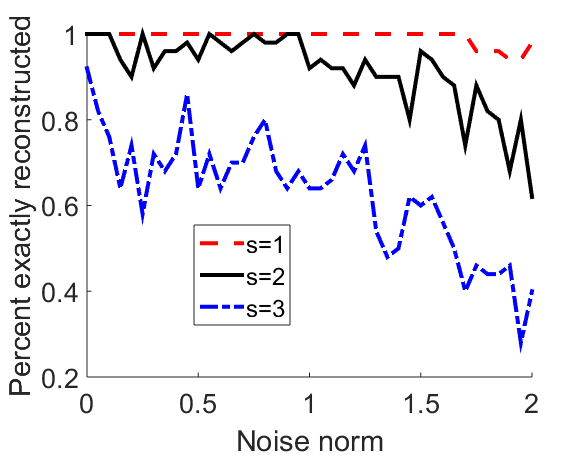}
\caption{Reconstruction results for the matrix described in Example \ref{ex4} using brute force closest vector problem (CVP) approach. Signals were generated with random support, and entries were integers uniform between -1 and 1 (left) or -5 and 5 (right). Gaussian noise was added to the measurements $y=Ax$ before reconstruction, at various levels. Note that our theory for this $3\times 6$ matrix only guarantees exact reconstruction of $s=1$-sparse signals and for noise at most $\sqrt{3}/2\approx 0.866$ (see Remark \ref{workitout}). Our experiments seem to show success slightly beyond this guarantee.}\label{therealfig2}
\end{figure}
\end{example}
\medskip

\begin{rem} \label{liouville} It is also possible to construct a $1 \times d$ algebraic matrix $A$ with $A\bx$ bounded away from $0$ for integer signals, but with the bound depending on $\bx$ and $d$. Specifically, if $K$ is a number field of degree $d$ with real algebraic numbers
$$1 = \alpha_1,\alpha_2,\dots,\alpha_d \in K$$
forming a $\que$-basis for $K$, let $A = (\alpha_1\ \dots\ \alpha_d)$ be the corresponding $1 \times d$ matrix. Let $1 \leq s \leq d$. Then for any $\bo \neq \bx \in \zed_s^d$,
$$|A \bx| > \frac{1}{C_s |\bx|^d},$$
where $C_s$ is an explicit constant depending on $s$, $d$, and $\alpha_1,\dots,\alpha_d$. This result follows from the argument of Section~3 of \cite{kron}, essentially constructing a class of badly approximable linear forms. More information on such linear forms can be found in~\cite{schmidt}, pp. 36--46. For some more recent information on this subject see, for instance~\cite{bennett}.
\end{rem}
\bigskip

\section{Reconstruction algorithm}
\label{recon}

Let us also say a few words about the reconstruction algorithm for our matrix construction in the situation when $s=m$. First recall the Closest Vector Problem (CVP) in some~$n$-dimensional Euclidean space~$\real^n$. This is an algorithmic lattice problem, which on the input takes a matrix $C \in \GL_n(\real)$ and a point $\bwy \in \real^n$ and on the output returns a point~$\bx$ in the lattice $\Lambda := C\zed^n$ such that
$$\|\bx-\bwy\| = \min \{ \|\bz - \bwy\| : \bz \in \Lambda \}.$$
Let the notation be as above with $s=m$, and $A$ be an $m \times d$ matrix with no $m$ column vectors linearly dependent, so that for all $\bx \in \zed_m^d$,
$$\|A\bx\| \geq \alpha$$
for some real $\alpha > 0$. Let $\J$ be the set of all cardinality $m$ subsets of $[d] = \{1,\dots,d\}$, then $|\J| = \binom{d}{m}$. For each $I \in \J$, let $A_I$ be the $m \times m$ submatrix of $A$ indexed by the elements of $I$ and let $\Lambda_I = A_I \zed^s$ be the corresponding lattice of rank~$m$ in~$\real^m$. Suppose now that $\bx \in \zed^d_m$, then $A\bx$ is a vector in some $\Lambda_I$. Let us write
$$\J = \{ I_1,\dots,I_t \},$$
where $t = \binom{d}{m}$. Given a CVP oracle, we can propose the following reconstruction algorithm for our problem.
\medskip

\noindent
{\bf Reconstruction Algorithm.} 
\begin{enumerate}

\item {\it Input:} A vector $\bwy = A\bx + \be$ for some $\bx \in \zed_m^d$ and error $\be \in \real^m$ with $\|\be\| < \alpha/2$.

\item {\it CVP:} Make $t$ calls to the CVP oracle in~$\real^m$ with the input $\Lambda_{I_j}$ and $\bwy$ for each $1 \leq j \leq t$; let 
$$\bz_1 \in \Lambda_{I_1},\dots,\bz_t \in \Lambda_{I_t}$$
be the vectors returned.

\item {\it Comparison:} Out of $\bz_1,\dots,\bz_t$, pick $\bz_i$ such that
$$\|\bz_i - \bwy\| < \alpha/2.$$
By our construction, there can be only one such vector.

\item {\it Matrix inverse:} Compute $(A_{I_i})^{-1}$.

\item {\it Reconstruction:} Take $\bx = (A_{I_i})^{-1} \bz_i$.

\end{enumerate}
\medskip

On the other hand, suppose we had an oracle for a reconstruction algorithm with the error bound~$\alpha$, call it RA. Given a point~$\bwy \in \real^m$, make a call to RA oracle, returning a vector~$\bx \in \zed_m^d$. Compute $\bz = A\bx$, then $\bz$ is in one of the lattices $\Lambda_{I_1},\dots,\Lambda_{I_t}$, and, assuming that~$\|\bz - \bwy\| < \alpha/2$, we have
$$\|\bz - \bwy\| = \min \left\{ \|\bu - \bwy\| : \bu \in \bigcup_{j=1}^t \Lambda_{I_j} \right\}.$$
Hence $\bz$ is a CVP solution for $\bwy$ in~$\bigcup_{j=1}^t \Lambda_{I_j}$. In other words, the problem of reconstructing the sparse signal from the image under such a matrix $A$ in~$\real^m$ has essentially the same computational complexity as CVP in~$\real^m$. It is known~\cite{CVP} that CVP in~$\real^m$ can be solved by a deterministic $O(2^{2m})$ time and $O(2^m)$ space algorithm, or by a randomized $2^{m+o(m)}$-time and space algorithm~\cite{CVP1}, which gives an idea of the complexity of our reconstruction algorithm.  Classical compressed sensing methods offer far more efficient complexity, but also require the sparsity level $s$ to be much less than $m$.  Our framework allows any $s\leq m$, which is a much taller order.
\bigskip

\section{Sparse geometry of numbers}
\label{sparse}

In this section we prove Theorem~\ref{main2}. We first recall Minkowski's Convex Body Theorem.

\begin{thm} \label{mink} Let $V$ be an $m$-dimensional subspace of $\real^d$, $1 \leq m \leq d$. Let $M$ be a convex $\bo$-symmetric body in $V$ and let $\Lambda \subset V$ be a lattice of full rank. Suppose that $\Vol_m(M) \geq 2^m \det \Lambda$. Then $M \cap \Lambda$ contains a nonzero point.
\end{thm}

\noindent
Let us also recall Vaaler's cube-slicing inequality (see Corollary to Theorem~1 of~\cite{vaaler}).

\begin{lem} \label{jeff} Let $\C^d(1)$ be a cube of sidelength 1 centered at the origin in $\real^d$, i.e.
$$\C^d(1) = \{ \bx \in \real^d : |x_i| \leq 1/2\ \forall\ 1 \leq i \leq d \}.$$
Let $V$ be an $m$-dimensional subspace of $\real^d$, $m \leq d$.  Then the $m$-dimensional volume of the section $\C_d(1) \cap V$~is 
$$\Vol_m(\C_d(1) \cap V) \geq 1.$$
\end{lem}

\noindent
We can use this lemma to prove a sparse version of Minkowski's Convex Body Theorem for parallelepipeds. As in Section~\ref{integer}, we write $[d] = \{1,\dots,d\}$, and whenever $A$ is a matrix with $d$ columns, we write $A_I$ for the submatrix of $A$ consisting of columns indexed by $I \subset [d]$.

\begin{prop} \label{mink1} Let $m \leq d$ be positive integers. Let $A \in \GL_d(\real)$, and let $P_A=A\C^d(1)$. Assume that for some $I \subset [d]$ with $|I| = m$,
\begin{equation}
\label{I_cond}
\sqrt{|\det (A_I^{\top} A_I)|} \geq 2^m.
\end{equation}
Then $P_A$ contains a nonzero point of $\zed_m^d$.
\end{prop}

\proof
Let $I \subset [d]$ with $|I|=m$ satisfy~\eqref{I_cond}, and  define
$$V_I = \{ \bx \in \real^d : x_j = 0\ \forall\ j \notin I \},$$
hence $V_I$ is an $m$-dimensional coordinate subspace of~$\real^d$. Then
$$A (\C^d(1) \cap A^{-1} V_I) = P_A \cap V_I,$$
and the $m$-dimensional volume of the section $P_A \cap V_I$ is
$$\Vol_m (P_A \cap V_I) = \sqrt{|\det (A_I^{\top} A_I)|}\ \Vol_m (\C^d(1) \cap A^{-1} V_I) \geq \sqrt{|\det (A_I^{\top} A_I)|} \geq 2^m,$$
by Lemma~\ref{jeff}. Now notice that $P_A  \cap V_I$ is a convex $\bo$-symmetric set in $V_I$. Let $\Lambda_I = V_I \cap \zed^d$ be the full-rank integer lattice in $V_I$, so $\det \Lambda_I = 1$. Hence, by Theorem~\ref{mink}, $P_A \cap V_I$ contains a nonzero point of $\Lambda_I$. Since $\Lambda_I \subset \zed_m^d$, this means that $P_A$ contains a nonzero point of $\zed_m^d$.
\endproof

\noindent
We can now prove a sparse version of Minkowski's Linear Forms Theorem.

\begin{thm} \label{mink_linear} Let $m \leq d$ be positive integers and $B \in \GL_d(\real)$. For each $1 \leq i \leq d$, let 
$$L_i(X_1,\dots,X_d) = \sum_{j=1}^d b_{ij} X_j$$
be the linear form with entries of the $i$-th row of $B$ for its coefficients. Let $c_1,\dots,c_d$ be positive real numbers such that for some $I = \{1 \leq j_1 < \dots < j_m \leq d\} \subset [d]$,
\begin{equation}
\label{I_cond-1}
c_{j_1} \cdots c_{j_m} \geq \left| \det \left( (B^{-1})_I^{\top} (B^{-1})_I \right) \right|^{-1/2}.
\end{equation}
Then there exists a nonzero point $\bx \in \zed_m^d$ such that
\begin{equation}
\label{Lix}
|L_{j_i} (\bx)| \leq c_{j_i}
\end{equation}
for each $1 \leq i \leq m$.
\end{thm}

\proof
Define a $d \times d$ diagonal matrix $D$ with diagonal entries $2c_1,\dots,2c_d$, and let $A = B^{-1}D$. Then
$$P_A = A\C^d(1) = \{ \bx \in \real^d : |L_i(\bx)| \leq c_i\ \forall\ 1 \leq i \leq d \}.$$
Let $I \subset [d]$ with $|I|=m$ satisfy~\eqref{I_cond-1}, and let $D(I)$ be the $m \times m$ diagonal matrix with diagonal entries $2c_{j_1},\dots,2c_{j_m}$. Notice that $A_I = (B^{-1})_I D(I)$, so 
$$A_I^{\top} A_I = D(I) \left( (B^{-1})_I^{\top} (B^{-1})_I \right) D(I).$$
Since $\det D(I) = 2^m c_{j_1} \cdots c_{j_m}$, \eqref{I_cond-1} implies that
$$\sqrt{|\det (A_I^{\top} A_I)|} = \left( \det D(I) \right) \sqrt{|\det \left( (B^{-1})_I^{\top} (B^{-1})_I \right)|} \geq 2^m.$$
Hence, by Proposition~\ref{mink1}, there exists a nonzero point $\bx \in P_A \cap \zed^d_m$, i.e. $\bx$ satisfies~\eqref{Lix}.
\endproof

\begin{cor} \label{mink_cor} Let $A$ be an $m \times d$ real matrix of rank $m \leq d$. Let $B \in \GL_d(\real)$ be a matrix whose first $m$ rows are the rows of $A$. Let $I= \{ 1,\dots,m \}$. Then there exists a nonzero point $\bx \in \zed_m^d$ such that
$$\|A\bx\| \leq \sqrt{m} \left| \det \left( (B^{-1})_I^{\top} (B^{-1})_I \right) \right|^{-1/2m}.$$
\end{cor}

\proof
For each $1 \leq i \leq m$ let
$$c_i = \left| \det \left( (B^{-1})_I^{\top} (B^{-1})_I \right) \right|^{-1/2m},$$
so $c_1 \cdots c_m = \left| \det \left( (B^{-1})_I^{\top} (B^{-1})_I \right) \right|^{-1/2}$. Then by Theorem~\ref{mink_linear}, for some point $\bx \in \zed_m^d$ inequality \eqref{Lix} holds with each $j_i=i$, and so
$$\| A\bx\|^2 = \sum_{i=1}^m L_i(\bx)^2 \leq m \left| \det \left( (B^{-1})_I^{\top} (B^{-1})_I \right) \right|^{-1/m}.$$
This completes the proof.
\endproof

\begin{rem} Let the notation be as in Corollary~\ref{mink_cor} above. If we write $\II_d$ and $\II_m$ for the $d \times d$ and $m \times m$ identity matrices, respectively, then $B B^{-1} = \II_d$, and so $A (B^{-1})_I = \II_m$. Using Cauchy-Binet formula, we have
$$1 = \det (A (B^{-1})_I) = \sum_{J \subset [d], |J| =m} \det (A_J) \det ((((B^{-1})_I)^{\top})_J),$$
while
$$\det \left( (B^{-1})_I^{\top} (B^{-1})_I \right) = \sum_{J \subset [d], |J| =m} \det ((((B^{-1})_I)^{\top})_J)^2.$$
\end{rem}
\smallskip

\proof[Proof of Theorem~\ref{main2}] Theorem~\ref{main2} now follows from Corollary~\ref{mink_cor}.
\endproof

We will now give a few examples of matrices with column vectors having equally large sup-norms for which the bound of Theorem~\ref{main2} is better than the naive bound~\eqref{naive}.

\begin{example} \label{ex_mink} We use notation of Theorem~\ref{main2}. Let $d=5$, $m=3$, and let 
$$A_1 = \begin{pmatrix} 15 & 15 & 4 & 13 & 15 \\ 2 & -1 & -15 & 2 & -13 \\ -13 & 2 & 1 & -15 & 4 \end{pmatrix},$$
then
$$A_1' = \begin{pmatrix} 3392/3905 & 23/355 & 3021/3905 \\ -1949/2130 & 3/710 & -1697/2130 \\ -6409/9372 & -19/284 & -5647/9372 \\ -6407/9372 & -17/284 & -6353/9372 \\ 13869/15620 & 1/1420 & 12047/15620 \end{pmatrix},$$
and so the bound of~\eqref{thm_2-4} is~$8.375...$, which is better than $25.980...$, the bound given by~\eqref{naive}.
\smallskip

Let $d=6$, $m=3$, and let
$$A_2 = \begin{pmatrix} 50000 & 20 & 40 & 3 & -50000 & 30 \\ -1 & -50000 & 20 & 40 & 4 & -50000 \\ -50000 & -1 & -50000 & -50000 & 20 & 40 \end{pmatrix},$$
then
$$A_2' = \begin{pmatrix} \ \\ \frac{3907968052500399551464}{269371733328769889476945} & \frac{782608564652549551187}{53874346665753977895389} & \frac{60146658957656226816}{4144180512750305991953} \\  \\ \frac{593868225682933391}{107748693331507955790778} & -\frac{780555233932686945}{53874346665753977895389} & \frac{22817064638544932}{4144180512750305991953} \\ \\ -\frac{31214424035248447494619}{2154973866630159115815560} & -\frac{3129803777458426735103}{215497386663015911581556} & -\frac{240538498302920543093}{16576722051001223967812} \\ \\ -\frac{167547159662404885}{9795335757409814162798} & \frac{14067072993904685}{4897667878704907081399} & -\frac{6446699511344665}{376743682977300544723} \\ \\ \frac{31195662429419099248023}{2154973866630159115815560} & \frac{3127927852920367185691}{215497386663015911581556} & \frac{240394120545528419881}{16576722051001223967812} \\ \\ -\frac{11261122160952583033}{1077486933315079557907780} & -\frac{1125763802241065145}{107748693331507955790778} & -\frac{86645392385558751}{8288361025500611983906} \\ \ \end{pmatrix},$$
and so the bound of~\eqref{thm_2-4} is~$7651.170...$, which is better than $86602.540...$, the bound given by~\eqref{naive}.
\smallskip

Let $d=8$, $m=4$, and let
$$A_3 = \begin{pmatrix} 6 & 13 & 13 & 11 & 6 & 12 & 11 & 10 \\ 7 & 12 & 6 & 13 & 7 & 11 & 11 & 9 \\ 8 & 11 & 12 & 9 & 12 & 12 & 12 & 11 \\ 13 & 10 & 7 & 8 & 13 & 13 & 13 & 13 \end{pmatrix},$$
then
$$A_3' = \begin{pmatrix} -736/1859 & 1865/5577 & 566/1859 & -661/1859 \\ 1990/1859 & 328/1859 & -3844/1859 & 1277/1859 \\ -1635/1859 & 646/1859 & 2495/1859 & -1577/1859 \\ -3015/1859 & 4273/5577 & 4021/1859 & -2584/1859 \\ 1499/1859 & -1654/5577 & -2350/1859 & 1273/1859 \\ 1228/169 & -2117/507 & -1523/169 & 1045/169 \\ -5605/1859 & 2243/1859 & 8286/1859 & -4218/1859 \\ -7461/1859 & 11917/5577 & 9390/1859 & -6289/1859 \end{pmatrix},$$
and so the bound of~\eqref{thm_2-4} is~$2.412...$, which is better than $26$, the bound given by~\eqref{naive}.

\end{example}

\bigskip

\bibliographystyle{myalpha}  
\bibliography{sparse-recovery-2}        

\begin{thebibliography}{CRTV05}

\bibitem[ADSD15]{CVP1}
D.~Aggarwal, D.~Dadush, and N.~Stephens-Davidowitz.
\newblock Solving the closest vector problem in $2^n$ time -- the discrete
  {G}aussian strikes again!
\newblock In {\em IEEE 56th Annual Symposium on Foundations of Computer Science
  -- FOCS 2015}, pages 563--582. IEEE Computer Soc., Los Alamitos, CA, 2015.

\bibitem[Alo99]{alon}
N.~Alon.
\newblock Combinatorial {N}ullstellensatz.
\newblock {\em Combin. Probab. Comput.}, 8(1-2):7--29, 1999.

\bibitem[AS08]{AS}
N.~Alon and J.~H. Spencer.
\newblock {\em The Probabilistic Method. Third eiditon}.
\newblock Wiley-Interscience, 2008.

\bibitem[BCV17]{cover}
M.~Balko, J.~Cibulka, and P.~Valtr.
\newblock Covering lattice points by subspaces and counting point-hyperplane
  incidences.
\newblock In B.~Aronov and M.~J. Katz, editors, {\em 33rd International
  Symposium on Computational Geometry (SoCG 2017)}, pages 12:1--12:16. Leibniz
  International Proceedings in Informatics Schloss Dagstuhl -- Leibniz-Zentrum
  fur Informatik, Dagstuhl Publishing, Germany, 2017.

\bibitem[BD09]{PaperIHT}
T.~Blumensath and M.~E. Davies.
\newblock Iterative hard thresholding for compressed sensing.
\newblock {\em Appl. Comput. Harmon. A.}, 27(3):265--274, 2009.

\bibitem[Ben96]{bennett}
M.~A. Bennett.
\newblock Simultaneous rational approximation to binomial functions.
\newblock {\em Trans. Amer. Math. Soc.}, 348(5):1717--1738, 1996.

\bibitem[BVW10]{random}
J.~Bourgain, V.~Vu, and P.~M. Wood.
\newblock On the singularity probability of discrete random matrices.
\newblock {\em J. Funct. Anal.}, 258(2):559--603, 2010.

\bibitem[CRT06]{RefWorks:47}
E.~J. Cand\`es, J.~Romberg, and T.~Tao.
\newblock Stable signal recovery from incomplete and inaccurate measurements.
\newblock {\em Commun. Pur. Appl. Math.}, 59(8):1207--1223, 2006.

\bibitem[CRTV05]{candes2005error}
E.~Candes, M.~Rudelson, T.~Tao, and R.~Vershynin.
\newblock Error correction via linear programming.
\newblock In {\em Foundations of Computer Science, 2005. FOCS 2005. 46th Annual
  IEEE Symposium on}, pages 668--681. IEEE, 2005.

\bibitem[CT05]{RefWorks:48}
E.~J. Cand\`es and T.~Tao.
\newblock Decoding by linear programming.
\newblock {\em IEEE T. Inform. Theory}, 51:4203--4215, 2005.

\bibitem[DR16]{davenport2016overview}
M.~A. Davenport and J.~Romberg.
\newblock An overview of low-rank matrix recovery from incomplete observations.
\newblock {\em IEEE Journal of Selected Topics in Signal Processing},
  10(4):608--622, 2016.

\bibitem[DT09]{RefWorks:288}
D.~L. Donoho and J.~Tanner.
\newblock Counting the faces of radomly-projected hypercubes and orthants, with
  applications.
\newblock {\em J.Amer.Math.Soc.}, 22(1):1--53, 2009.

\bibitem[EK12]{RefWorks:373}
Y.~C. Eldar and G.~Kutyniok.
\newblock {\em Compressed sensing: theory and applications}.
\newblock Cambridge University Press, 2012.

\bibitem[FK17]{flinth2017promp}
A.~Flinth and G.~Kutyniok.
\newblock Promp: A sparse recovery approach to lattice-valued signals.
\newblock {\em Applied and Computational Harmonic Analysis}, 2017.

\bibitem[FM17]{kron}
L.~Fukshansky and N.~G. Moshchevitin.
\newblock On an effective variation of {K}ronecker's approximation theorem
  avoiding algebraic sets.
\newblock {\em preprint}, 2017.

\bibitem[FR13]{RefWorks:45}
S.~Foucart and H.~Rauhut.
\newblock {\em A Mathematical Introduction to Compressive Sensing}.
\newblock Springer, 2013.
\newblock In press.

\bibitem[MR11]{mangasarian2011probability}
O.~L. Mangasarian and B.~Recht.
\newblock Probability of unique integer solution to a system of linear
  equations.
\newblock {\em European Journal of Operational Research}, 214(1):27--30, 2011.

\bibitem[MV13]{CVP}
D.~Micciancio and P.~Voulgaris.
\newblock A deterministic single exponential time algorithm for most lattice
  problems based on {V}oronoi cell computations.
\newblock {\em SIAM J. Comput.}, 42(3):1364--1391, 2013.

\bibitem[NT09]{NeedeT_CoSaMP}
D.~Needell and J.~Tropp.
\newblock {CoSaMP}: {I}terative signal recovery from incomplete and inaccurate
  samples.
\newblock {\em Appl. Comput. Harmon. A.}, 26(3):301--321, 2009.

\bibitem[RHE14]{rossi2014spatial}
M.~Rossi, A.~M. Haimovich, and Y.~C. Eldar.
\newblock Spatial compressive sensing for mimo radar.
\newblock {\em IEEE Transactions on Signal Processing}, 62(2):419--430, 2014.

\bibitem[Sch80]{schmidt}
W.~M. Schmidt.
\newblock {\em Diophantine Approximation}.
\newblock Lecture Notes in Mathematics, 785. Springer, Berlin, 1980.

\bibitem[Sto10]{stojnic2010recovery}
M.~Stojnic.
\newblock Recovery thresholds for ℓ 1 optimization in binary compressed
  sensing.
\newblock In {\em Information Theory Proceedings (ISIT), 2010 IEEE
  International Symposium on}, pages 1593--1597. IEEE, 2010.

\bibitem[TG07]{Paper9}
J.~A. Tropp and A.~C. Gilbert.
\newblock Signal recovery from random measurements via {O}rthogonal {M}atching
  {P}ursuit.
\newblock {\em IEEE T. Inform. Theory}, 53(12):4655--4666, 2007.

\bibitem[TLL09]{tian2009detection}
Z.~Tian, G.~Leus, and V.~Lottici.
\newblock Detection of sparse signals under finite-alphabet constraints.
\newblock In {\em IEEE International Conference on Acoustics, Speech and Signal
  Processing (ICASSP)}, pages 2349--2352. IEEE, 2009.

\bibitem[Vaa79]{vaaler}
J.~D. Vaaler.
\newblock A geometric inequality with applications to linear forms.
\newblock {\em Pacific J. Math.}, 83(2):543--553, 1979.

\bibitem[ZG11]{zhu2011exploiting}
H.~Zhu and G.~B. Giannakis.
\newblock Exploiting sparse user activity in multiuser detection.
\newblock {\em IEEE Transactions on Communications}, 59(2):454--465, 2011.

\bibitem[Zha11]{RefWorks:150}
T.~Zhang.
\newblock Sparse recovery with orthogonal matching pursuit under {RIP}.
\newblock {\em IEEE T. Inform.Theory,}, 57(9):6215--6221, 2011.

\end{thebibliography}

\end{document}